\documentclass[11pt,a4paper,twoside]{amsart}
\usepackage{amsmath}
\usepackage{amssymb}
\usepackage{latexsym}

\newcommand{\co}{\colon\thinspace}    
\newcommand{\fnote}[1]{\footnote{\small sharp1}}
\newcommand{\inv}{^{-1}}              

\newcommand{\N}{{\mathbb N}}
\newcommand{\Z}{{\mathbb Z}}
\newcommand{\R}{{\mathbb R}}

\newcommand{\T}{{\mathbb T}}

\newtheorem{theorem}{Theorem}[section]

\newtheorem{lemma}[theorem]{Lemma}
 
\newtheorem{conjecture}[theorem]{Conjecture} 

\title{Two remarks about Ma\~n\'e's conjecture}
\author{Daniel Massart}
\date{\today}
\begin{document}


\maketitle

\section{Introduction}
In this note we consider an autonomous Tonelli Lagrangian $L$ on a closed manifold $M$, that is, a $C^2$ function $L \co TM \longrightarrow \R$ such that $L$ is fiberwise strictly convex and superlinear. Then the Euler-Lagrange equation associated with $L$ defines a complete flow $\phi_t$ on $TM$. 
Define $\mathcal{M}_{inv}$
to be the set of $\Phi_t$-invariant, compactly supported, Borel probability measures on $TM$.
Mather showed that the function (called action of the Lagrangian on measures)
	\[
	\begin{array}{rcl}
\mathcal{M}_{inv} & \longrightarrow & \R \\
\mu & \longmapsto & \int_{TM}	L d\mu
\end{array}
\]
is well defined and has a minimum.  A measure achieving this minimum is called $L$-minimizing. The union, in $TM$, of the support of all minimizing measures is called Mather set of $L$, and denoted $\mathcal{M}(L)$. It is  compact and $\phi_t$-invariant.  See \cite{Mather91} and \cite{Fathi_bouquin} for more background. 

Observe that if $f$ is a $C^2$ function on $M$, $L+f$ is  also a  Tonelli Lagrangian. Adding a function to a Lagrangian is called perturbing the Lagrangian by a potential. Following Ma\~n\'e we say a property holds for a generic Lagrangian if, given any Lagrangian, the property holds for a generic perturbation by a potential. 
Ma\~n\'e conjectured a generic description of the minimizing measures : 
\begin{conjecture} [\cite{Mane97}] \label{forte}
Let
\begin{itemize}
	\item $M$ be a closed manifold
	\item $L$ be an autonomous Tonelli Lagrangian on $TM$
	
	\item $\mathcal{O}_1(L)$ be the set of $f$ in $C^{\infty}(M)$ such that the Mather set of $L+f$ consists of one periodic orbit.
	\end{itemize}
Then  the set $\mathcal{O}_1(L)$ is residual in $C^{\infty}(M)$.
\end{conjecture}
In other words, for a generic Lagrangian, there exits a unique minimizing measure, and it is supported by a periodic orbit. A similar conjecture can be made replacing $C^{\infty}(M)$ by $C^{k}(M)$ for any $k \geq 2$. 

Many more interesting invariant sets can be obtained by minimization than just the Mather set. If $\omega$ is a closed one-form on $M$, then $L-\omega$ is a Tonelli Lagrangian, and it has the same Euler-Lagrange flow as $L$. Its Mather set, however, is different in general. The Mather set of $L-\omega$ only depends on the cohomology class $c$ of $\omega$, we denote it $\mathcal{M}(L,c)$. It is often interesting to obtain information simultaneously on the Mather sets $\mathcal{M}(L,c)$ for a large set of cohomology classes. Thus Ma\~n\'e proposed the 
\begin{conjecture}[\cite{Mane96}]\label{faible}
 If $L$ is a Tonelli Lagrangian  on a manifold $M$, there exists  a residual subset $\mathcal{O}_2(L)$ of  $C^{\infty}(M)$, such that for any $f$ in $\mathcal{O}_2(L)$, there exists an open and dense subset $U(L,f)$ of $H^1 (M,\R)$ such that, for any $c$ in $U(L,f)$, the Mather set of $(L,c)$ consists of one periodic orbit.
\end{conjecture}
Intuitively Conjecture \ref{faible} is weaker than Conjecture \ref{forte} because we allow a larger set of perturbations (potentials and closed one-forms instead of just potentials). However the requirement of an open dense set in Conjecture \ref{faible} makes it far from obvious. In section \ref{section2} we prove that Conjecture \ref{forte} contains Conjecture \ref{faible}, using recent tools from Fathi's weak KAM theory, the most prominent of which is the Aubry set $\mathcal{A}(L)$. All we need to know about the Aubry set is that
\begin{itemize}
	\item it consists of the Mather set, and (possibly) orbits homoclinic to the Mather set (see \cite{Fathi_bouquin})
	\item when there is only one minimizing measure, the Aubry set is upper semi-continuous as a function of the Lagrangian, that is, for any neighborhood $V$ of $\mathcal{A}(L)$ in $TM$, there exists a neighborhood $\mathcal{U}$ of $L$ in the $C^2$ compact-open topology, such that for any $L_1$ in $\mathcal{U}$, we have $\mathcal{A}(L_1)\subset V$ (see \cite{Bernard_Conley}).
\end{itemize}
We first prove that Conjecture \ref{forte} is equivalent to the apparently stronger 
\begin{conjecture}\label{forte_Aubry}
Let
\begin{itemize}
	\item $M$ be a closed manifold
	\item $L$ be an autonomous Tonelli Lagrangian on $TM$
	\item $\mathcal{O}_3(L)$ be the set of $f$ in $C^{\infty}(M)$ such that the Aubry set of $L+f$ consists of one, hyperbolic periodic orbit.
\end{itemize}
Then the set $\mathcal{O}_3(L)$ is residual in $C^{\infty}(M)$.
\end{conjecture}
Then we prove that Conjecture \ref{forte_Aubry} contains the following, which obviously contains Conjecture \ref{faible} :
\begin{conjecture}\label{faible_Aubry}
 If $L$ is a Tonelli Lagrangian  on a manifold $M$, there exists  a residual subset $\mathcal{O}_4(L)$ of  $C^{\infty}(M)$, such that for any $f$ in $\mathcal{O}_2(L)$, there exists an open and dense subset $U(L,f)$ of $H^1 (M,\R)$ such that, for any $c$ in $U(L,f)$, the Aubry set of $(L+f,c)$ consists of one, hyperbolic periodic orbit.
\end{conjecture}
Conjecture \ref{faible_Aubry} is proved, in the case where the dimension of $M$ is two, in \cite{lowdim} (after a sketch of a proof appeared in \cite{ijm}). The analogous statement for Lagrangians which depend periodically on time is proved, in the case where the dimension of $M$ is one, in \cite{Osvaldo}. 

Conjecture \ref{forte} may be seen as an Aubry-Mather version of the Closing Lemma. 
This suggests that it should be true in the $C^2$ topology on Lagrangians and false in the $C^k$ topology for $k > 2$. If we want to prove the $C^k$ version of Conjecture \ref{forte}, and we are lucky enough to have a sequence of periodic orbits $\gamma_n$ which approximate our Mather set, then the first idea that comes to mind is to perturb $L$ by a non-negative potential $f_n$ which vanishes only on $\gamma_n$. Then $\gamma_n$ is still an orbit of $L+f_n$. If we can find $f_n$ big enough for $\gamma_n$ to be $L+f_n$-minimizing, but small enough for the $C^k$-norm of $f_n$ to converge to zero, then we are done. In Section \ref{section3} we prove that this naive approach doesn't work in the $C^k$-topology, for $k \geq 4$. Specifically, we give an example of a Lagrangian $L$ on the two-torus, such that for any periodic orbit $\gamma$ of $L$, and any $C^4$ function $f$ on the two-torus, if $\gamma$ is $L+f$-minimizing, then the $C^4$ norm of $f$ is bounded below by a constant which only depends on $L$. 

\textbf{Acknowledgements} This work was partially supported by the ANR project ''Hamilton-Jacobi et th\'eorie KAM faible''.
\section{}\label{section2}
\begin{lemma}\label{lemme1}
Let
\begin{itemize}
	\item $M$ be a closed manifold
	\item $L$ be an autonomous Tonelli Lagrangian on $TM$
	\item $\mathcal{O}_3(L)$ be the set of $f$ in $C^{\infty}(M)$ such that the Aubry set of $L+f$ consists of one, hyperbolic periodic orbit
	\item $\mathcal{O}_1(L)$ be the set of $f$ in $C^{\infty}(M)$ such that the Mather set of $L+f$ consists of one periodic orbit.
\end{itemize}
Then $\mathcal{O}_3(L)$ is open and dense in $\mathcal{O}_1(L)$.
\end{lemma}
\proof
We first prove that $\mathcal{O}_3(L)$ is open  in $\mathcal{O}_1(L)$.
Take $f \in \mathcal{O}_3(L)$. Replacing $L$ with $L+f$, we may assume $f=0$. Let $\gamma$ be the hyperbolic periodic orbit which comprises $\mathcal{A}(L+f)$.
By a classical property of hyperbolic periodic orbits,  there exists a neighborhood $\mathcal{U}_1$ of the zero function in $C^{\infty}(M)$, and a neighborhood $V$ of $\gamma$ in $TM$ such that for any $f \in \mathcal{U}_1$, for any energy level $E$ of $L$, the only invariant set of the Euler-Lagrange flow of $L$ contained in $E \cap V$, if any, is a hyperbolic periodic orbit homotopic to $\gamma$.

Since $\mathcal{A}(L)$ is a  periodic orbit, the quotient Aubry set $A$ has but one element. Thus by \cite{Bernard_Conley}, there exists a neighborhood $\mathcal{U}_2$ of the zero function in $C^{\infty}(M)$, such that for all $f$ in $\mathcal{U}_2$, we have $\mathcal{A}(L+f)\subset V$. Therefore, for any $ f \in \mathcal{U}_1 \cap \mathcal{U}_2$, the Aubry set $\mathcal{A}(L+f)$ consists of one,  hyperbolic periodic orbit. 

Now let us prove that $\mathcal{O}_3(L)$ is dense  in $\mathcal{O}_1(L)$.

Take $f \in \mathcal{O}_1(L)$. Replacing $L$ with $L+f$, we may assume $f=0$. Let $\gamma$ be the periodic orbit which comprises $\mathcal{M}(L)$.
Now let us take a smooth function $g$ on $M$ such that $g$ vanishes on the projection to $M$ of $\gamma$ (which we again denote $\gamma$ for simplicity), and $\forall x \in M,\  g(x) \geq d(x,\gamma)^2$, where the distance is meant with respect to some Riemannian metric on $M$. Let $\lambda$ be any positive number. We will show that $\lambda g \in \mathcal{O}(h)$, which proves that $\mathcal{O}(h)$ is dense  in $\mathcal{O}_1(L)$. Observe that  $\gamma$ is a minimizing hyperbolic periodic orbit of the Euler-Lagrange flow of $L+\lambda g$ (see \cite{CI99}). 
Furthermore , $\alpha_{L+\lambda g}(0) =\alpha_L(0)$, where $\alpha_L(0)$ is Ma\~n\'e's critical value for the Lagrangian $L$. 

Adding a constant to $L$ if necessary, we assume $\alpha_L(0)=0$. Recall that the Aubry set is the union of the Mather set and orbits homoclinic to the Mather set. Therefore, to prove that $\lambda g \in \mathcal{O}_0 (h)$, it suffices to prove that the Aubry set 
$\mathcal{A}\left(L+\lambda g \right)$ does not contain any orbit homoclinic to  $\gamma$.

Assume $\delta \co \R \longrightarrow M$ is an extremal of $L+\lambda g$, homoclinic to  $\gamma$. Since $g(\delta(t)) >0$ for all $t$, there exists $C>0$ such that
	\[ \int^{+\infty}_{-\infty} g(\delta (t))dt  \geq 2C.
\]
Let $u$ be a weak KAM solution for $L$. We have, for any $s,t \in \R$,
 remembering that $\alpha_L(0)=0$, 
	\[ \int^{t}_{s} L(\delta (t),\dot{\delta}(t))dt \geq u\left(\delta(t)\right)-u\left(\delta(s)\right).
\]
Since is homoclinic to  $\gamma$ there exist two sequences $t_n$ and $s_n$ that converge to $+\infty$, such that $\delta(t_n)$ and $\delta (-s_n)$ converge to the same point $x$ on $\gamma$, so for $n$  large enough
	\[ \int^{t_n}_{-s_n} L(\delta (t),\dot{\delta}(t))dt \geq -C.
\]
 Therefore, for $n$ large enough,  
	\[\int^{t_n}_{-s_n} \left(L+\lambda g \right)(\delta (t),\dot{\delta}(t))dt > C.
\]
On the other hand, since $\alpha_{L+\lambda g}(0) =0$, if $\delta$ were contained in the projected Aubry set of $L +\lambda g $, we would have, denoting $u_{\lambda}$ a weak KAM solution for $L +\lambda g $,
	\[\int^{t_n}_{-s_n} \left(L+\lambda g \right)(\delta (t),\dot{\delta}(t))dt = u_{\lambda}\left(\delta(t)\right)-u_{\lambda}\left(\delta(s)\right)
\]
which converges to zero because $\delta(t_n)$ and $\delta (-s_n)$ converge to the same point $x$ on $\gamma$.
 Therefore the Aubry set of $L +\lambda g $ consists of  $\gamma$ alone, which proves that $\lambda g \in \mathcal{O}_3(L)$, and the Lemma. \qed 

Therefore $\mathcal{O}(L)$ is residual in $C^{\infty}(M)$ if and only if $\mathcal{O}_1(L)$ is. Therefore, Conjecture \ref{forte} is equivalent to Conjecture \ref{forte_Aubry}. Now we show that Conjecture \ref{forte_Aubry} contains Conjecture \ref{faible_Aubry}, which obviously contains Conjecture \ref{faible}. 

Assume Conjecture \ref{forte_Aubry} is true.
Let $L$ be an autonomous Tonelli Lagrangian on a manifold $M$. Let $c_i, i \in \N$ be a countable dense subset of $H^1 (M,\R)$. Take, for every $i \in \N$, a closed one-form $\omega_i$ with cohomology $c_i$. Since Conjecture \ref{forte_Aubry} is true for every Lagrangian $L-\omega_i$,  for every $i \in \N$, there exists a residual 
subset $\mathcal{O}_i$ of $C^{\infty}(M)$ such that for every $f$ in $\mathcal{O}_i$, $\mathcal{A}(L+f,c_i)$ consists of one hyperbolic periodic orbit. 
Then the intersection, over $i \in \N$, of $\mathcal{O}_i$ is a residual subset $\mathcal{O}$ of $C^{\infty}(M)$. For every $f$ in $\mathcal{O}$, for every $i \in \N$, $\mathcal{A}(L+f,c_i)$ consists of one hyperbolic periodic orbit $\gamma_i$. As in the proof of Lemma \ref{lemme1}, there exists a neighborhood $V_i$ of $c_i$ in $H^1(M,\R)$, such that for any $c$ in $V_i$, the Aubry set $\mathcal{A}(L+f,c)$ consists of one hyperbolic periodic orbit homotopic to $\gamma_i$. The union, over $i \in \N$, of the $V_i$, is an open and dense subset $V$ of  $H^1(M,\R)$, and for any $c$ in $V$, $\mathcal{A}(L+f,c)$ consists of one hyperbolic periodic orbit, so Conjecture \ref{faible_Aubry} is true.

\section{An example}\label{section3}
Let
\begin{itemize}
	\item $r$ be a quadratic irrational number, for instance $\sqrt{2}$
	\item $p_0$ and $q_0$ be real numbers such that $p^{2}_{0}+q^{2}_{0}=1$ and $p_0 / q_0 = r$
	\item $\T^2$ be $\R^2 / \Z^2$, endowed with canonical coordinates $(x,y)$
	\item $L$ be the Lagrangian on $T\T^2$ defined by 
	\[ L(x,y,u,v) := \frac{u^2 + v^2}{2}- \left(p_0 u + q_0 v\right)
\]
where $(u,v)$ are the tangent coordinates to $(x,y)$.
\end{itemize}
Assume that for some function $f$ on $\T^2$, $L+f$ has a minimizing periodic orbit $\gamma$, and furthermore, $\gamma$ is an orbit of $L$, that is, it has the form $t \longmapsto (pt,qt)$ for some real numbers $p$ and $q$. Then, if $T$ is the smallest period of $\gamma$, $(pT,qT) \in \Z^2$ and $pT$, $qT$ are mutually prime. Consider the map
	\[
	\begin{array}{rcl}
	F \co \R & \longrightarrow & \R \\
	\lambda & \longmapsto & \frac{1}{T}\int^{T}_{0} f(pt, qt + \lambda) dt.
	\end{array}
\]
Observe that $F$ is $1$-periodic. We now prove that $F$ is $(pT)\inv$-periodic. Indeed, take $r,s$ in $\Z$ such that $pTr-qTs=1$. Then for any $t$,
	\begin{eqnarray*} (pt, qt + \frac{1}{pT})& = & (pt, qt + r -\frac{qs}{p})\\
	& = & \left(pt, q(t -\frac{s}{p}) \right) \  \mbox{mod}\Z^2 \\
	&=& \left(p(t -\frac{s}{p})+s, q(t -\frac{s}{p}) \right) \  \mbox{mod}\Z^2 \\
	&=& \left(p(t -\frac{s}{p}), q(t -\frac{s}{p}) \right) \  \mbox{mod}\Z^2 \\
\end{eqnarray*}
so
	\begin{eqnarray*} 
	F\left(\lambda + \frac{1}{pT}\right)&=& \frac{1}{T}\int^{T}_{0} f(pt, qt + \frac{1}{pT}+ \lambda) dt \\
&=& \frac{1}{T}\int^{T}_{0} f(p(t -\frac{s}{p}), q(t -\frac{s}{p})+ \lambda) dt = F\left(\lambda \right)
\end{eqnarray*}
using the change of variable $t \mapsto t -s/p$ (and the fact that $F$ is $1$-periodic). Now we prove  that
	\[ \int^{1}_{0}F(\lambda) d\lambda = \int_{\T^2}f d\mbox{leb}
\]
where $\mbox{leb}$ denotes the standard Lebesgue measure on $\T^2=\R^2 / \Z^2$. Indeed, let $\nu$ be  the measure on $\T^2$ defined by 
	\[ \int g (x,y) d\nu (x,y) := \int^{1}_{0} d\lambda \left\{ \frac{1}{T}\int^{T}_{0} g(pt,qt + \lambda) dt \right\}
\]
for any continuous function $g$ on $\T^2$.
We want to prove that $\nu $ is actually $\mbox{leb}$. First let us show that $\nu$ is invariant under translations. Let $(u,v)$ be any vector in $\R^2$. We have
\begin{eqnarray*} 
\int g(x+u,y+v) d\nu (x,y) &=& \int^{1}_{0} d\lambda \left\{ \frac{1}{T}\int^{T}_{0} g(pt+u,qt + \lambda +v ) dt \right\}  \\
&=& \int^{1}_{0} d\lambda \left\{ \frac{1}{T}\int^{T}_{0} g(p(t+\frac{u}{p}),q(t+\frac{u}{p}) + \lambda - \frac{uq}{p}) dt \right\}  \\
&=& \int^{1}_{0} d\lambda \left\{ \frac{1}{T}\int^{T}_{0} g(pt,qt + \lambda - \frac{uq}{p}) dt \right\}  \\
&=& \int^{1}_{0} d\lambda \left\{ \frac{1}{T}\int^{T}_{0} g(pt,qt ) dt \right\}  = \int g (x,y) d\nu (x,y)\\
\end{eqnarray*}
where we have used, in succession, the changes of variables $t \mapsto t +u/p$ and $\lambda \mapsto \lambda - uq/p$. So $\nu$ is invariant under translations. Furthermore
$\int 1 d\nu = 1$ so $\nu $ is actually $\mbox{leb}$. 

Now let us use the fact that $\gamma$ is $L+f$-minimizing. Let $\mu$ be the probability measure equidistributed along $\gamma$. We have
\begin{eqnarray*}
\int \left(L+f\right) d\mu &=& \frac{1}{T}\int^{T}_{0} dt \left\{ \frac{p^2 + q^2}{2}- \left(p_0 p + q_0 q\right) +f(pt,qt)  \right\} \\
&=&  \frac{p^2 + q^2}{2}- \left(p_0 p + q_0 q\right) +F(0).
\end{eqnarray*}
Let $\mu_0$ be the measure on $T\T^2$ defined by 
	\[ \int g (x,y,u,v) d\mu_0 (x,y,u,v) := \int^{1}_{x=0}dx \; \int^{1}_{y=0} g(x,y,p_0,q_0) dy
\]
for any continuous function $g$ on $T\T^2$. Observe that  the measure $\mu_0$ is $L$-minimizing. In particular it is closed (see \cite{FS}, Theorem 1.6).  We have 
	\[ \int_{T\T^2} L d\mu_0 = -\frac{1}{2} \mbox{ and } \int_{T\T^2} f d\mu_0 = \int_{\T^2} f d\mbox{leb}= \int^{1}_{0}F(\lambda) d\lambda 
\]
where we have implicitely extended $f$ to a function on $T\T^2$ by setting $f(x,y,u,v):=f(x,y)$ for any $u$ and $v$. 
Since $\mu$ is $L+f$-minimizing, and $\mu_0$ is closed, we have (see \cite{FS}, Theorem 1.6) 
	\[ \int_{T\T^2} (L+f) d\mu_0 \geq \int_{T\T^2} (L+f) d\mu
\]
that is, 
\begin{eqnarray*}
\int^{1}_{0}F(\lambda) d\lambda -F(0) &  \geq & \int_{T\T^2} L d\mu -\int_{T\T^2} L d\mu_0 \\
&=& \frac{p^2 + q^2}{2}- \left(p_0 p + q_0 q\right) +\frac{1}{2} \\
&=& \frac{1}{2}(p-p_0)^2 + \frac{1}{2}(q-q_0)^2 .
\end{eqnarray*}
Now let us use the fact that $r=p_0 /q_0$ is quadratic, and $pT/qT = p/q$ is rational, so there exists a constant $C_0$ such that 
	\[\left| \frac{p_0}{q_0}-\frac{p}{q}\right| \geq \frac{C_0}{(pT)^2}.
\]
Hence, setting $C := C^{2}_{0}/2$,
	\[ \frac{1}{2}(p-p_0)^2 + \frac{1}{2}(q-q_0)^2 \geq \frac{C}{(pT)^4}
\]
whence
	\[\int^{1}_{0}F(\lambda) d\lambda -F(0) \geq \frac{C}{(pT)^4}.
\]
Therefore, since $F$ is $(pT)\inv$-periodic, there exists a $\lambda_0 \in \left[0, (pT)\inv \right]$ such that 
	\[ F(\lambda_0) -F(0) \geq \frac{C}{(pT)^4}.
\]
Hence there exists a $\lambda_1 \in \left[0, (pT)\inv \right]$ such that $\left|F'(\lambda_1)\right| \geq C(pT)^{-3}$. On the other hand, since $F$ is $(pT)\inv$-periodic, there exists a $\lambda_2  \in \left[0, (pT)\inv \right]$ such that $F'(\lambda_2)=0$. Thus $ \left| F'(\lambda_1)-F'(\lambda_2)\right|\geq C(pT)^{-3}$, so there exists a $\lambda_3 \in \left[0, (pT)\inv \right]$ such that $\left|F''(\lambda_3)\right| \geq C(pT)^{-2}$. Iterating this process we show there exists a $\lambda \in \left[0, (pT)\inv \right]$ such that $\left|F^{(4)}(\lambda)\right| \geq C$, that is, 
	\[\left|\frac{1}{T}\int^{T}_{0} f^{(4)}(pt,qt+\lambda ) dt  \right|\geq C.
\]
In particular the $C^4$-norm of $f$ is bounded below by $C$.

\end{document}